\documentclass[12pt]{article}
\usepackage[latin1]{inputenc}
\usepackage{amsmath}
\usepackage{amssymb}
\usepackage{graphicx}
\usepackage{graphics}
\usepackage{amsthm}
\usepackage{amscd}
\usepackage{color}
\usepackage{amsfonts}
\usepackage{fancyhdr}

\newtheorem{theorem}{Theorem}[section]
\newtheorem{mainthm}{Theorem}
\newtheorem{Co}[mainthm]{Corollary}
\newtheorem*{theorem*}{Theorem}
\newtheorem{corollary}[theorem]{Corollary}
\newtheorem{proposition}[theorem]{Proposition}
\newtheorem{lemma}[theorem]{Lemma}
\newtheorem{Remark}[theorem]{Remark}
\newtheorem*{definition*}{Definition}

\newtheorem{claim}[theorem]{Claim}

\newtheorem*{Question*}{Question}
\newtheorem{Question}{Question}

\def\N{\mathbb{N}}

\def\e{\varepsilon}

\def\norm #1{\Vert \,#1\, \Vert\,}
\makeatletter

\newcommand{\Rmnum}[1]{\expandafter\@slowromancap\romannumeral #1@}
\makeatother

\def\ud{\mathrm{d}}

\def\diff {\operatorname{Diff}}
\def\dim{\operatorname{dim}}
\def\Cov{\operatorname{Cov}}

\def\La{\Lambda}

   \def\cM{\mathcal{M}} 
    
\def\cC{\mathcal{C}}   \def\cO{\mathcal{O}} \def\cU{\mathcal{U}}
    \def\cV{\mathcal{V}}
    
\def\cF{\mathcal{F}}

\def\diam{\operatorname{diam}}

\begin{document}
\vspace{-2cm}
\title {Non-hyperbolic ergodic measures and horseshoes in partially hyperbolic homoclinic classes}
\author{Dawei YANG \and Jinhua ZHANG\footnote{This work was done when J. Zhang visited Soochow University in July 2017. J. Zhang would like to thank Soochow University for hospitality. 
	D. Yang  was partially supported by NSFC 11671288 and NSFC 11790274. J. Zhang was partially supported by the ERC project 692925 NUHGD. J. Zhang is the corresponding author.}}

\maketitle

\begin{abstract}
	We study a rich family of robustly non-hyperbolic transitive  diffeomorphisms and we show that  each   ergodic measure is approached by hyperbolic sets  in weak$*$-topology and in entropy. For hyperbolic ergodic measures, it is a classical result of A. Katok. The novelty here is to deal with non-hyperbolic ergodic measures.

\end{abstract}
\begin{flushleft}
\end{flushleft}
\begin{flushleft}
\end{flushleft}

\section{Introduction}

Given a diffeomorphism $f$ on   a smooth Riemannian manifold $M$, one says that $f$ is \emph{partially hyperbolic}, if there exist  an invariant splitting $TM=E^s\oplus E^c\oplus E^u$ and  a metric $\|\cdot\|$ such that for each $x\in M$, one has
$$\norm{Df|_{E^s(x)}}<\min\{1,m(Df|_{E^c(x)})\}\leq\max\{1,\norm{Df|_{E^c(x)}}\}<m(Df|_{E^u(x)}).$$

In this paper, we consider the 
  set  $\cU(M)$ of  all  $C^1$ partially hyperbolic diffeomorphisms on $M$  satisfying that for each $f\in\mathcal{U}(M)$, one has:
\begin{itemize}
\item $f$ is partially hyperbolic with 1-dimensional center bundle;
\item $f$ has hyperbolic periodic points of different indices;
\item the strong stable and unstable foliations are robustly minimal.
\end{itemize}
By definition, $\cU(M)$ is an open set.

There are two types of  classical examples of partially hyperbolic diffeomorphisms with two minimal strong foliations  due to ~\cite{BDU}:
\begin{itemize} 
	\item   the time one map of a transitive Anosov flow admits  perturbations with robustly minimal strong foliations. 
	\item  the skew-product of  the rotations on circle over  an Anosov diffeomorphism on torus admits  perturbations with robustly minimal strong foliations. 
	\end{itemize}
 Very recently, some new examples of partially hyperbolic diffeomorphisms with two minimal strong foliations are constructed  and such examples are not isotopic to the classical ones listed above,  see the anomalous examples in  ~\cite{BGP,BGHP}.  

 Before stating our result, let us recall some notions.  Given a compact $f$-invariant set $K\subset M$, one denotes by $\cM_{inv}(K,f)$ the set of $f$-invariant measures supported on $K$.
A compact invariant set $K$ is \emph{hyperbolic} if one has the splitting $T_KM=E^s\oplus E^u$ such that $E^s$ is uniformly contracting and $E^u$ is uniformly expanding; and the \emph{index of a hyperbolic set} is defined as the dimension of its stable bundle. 
Recall that  an ergodic measure $\mu$ is \emph{hyperbolic}, if all the Lyapunov exponents of $\mu$ are non-zeo;  and we say $\mu$ is \emph{non-hyperbolic}, if $\mu$ has vanishing Lyapunov exponents. In our setting, we succeed in approaching   ergodic measures by hyperbolic sets  in weak$*$-topology as well as in entropy.

\begin{mainthm}\label{thmA} 
	There exists a $C^1$ open and dense subset $\mathcal{V}(M)$ of $\mathcal{U}(M)$ such that for any $f\in\mathcal{V}(M)$,  each  $f$-ergodic measure $\mu$ is approached by hyperbolic sets  \emph{in weak$*$-topology and in entropy},
	that is,  for any $\e>0$, there exists a hyperbolic set  $\Lambda_{\varepsilon}$   such that
	\begin{itemize}
		\item[--] 
		$|h_{\mu}(f)-h_{top}(f|_{\Lambda_{\varepsilon}})|<\varepsilon;$
		\item[--] for any  $\nu\in\cM_{inv}(\Lambda_{\varepsilon},f)$, one has $\ud(\mu,\nu)<\varepsilon$. 
\end{itemize}
Moreover,  if $\mu$ is non-hyperbolic, $\mu$ can be  approxiamted by hyperbolic sets of different indices  in weak$*$-topology and in entropy, i.e. the index of $\La_\e$ can be chosen to be $\dim(E^s)$ or $\dim(E^s)+1$.
\end{mainthm}

 In  $C^1$-setting, if the support of  a hyperbolic ergodic measure $\nu$   admits a  dominated  splitting which coincides with the hyperbolic splitting of the measure $\mu$, then  $\mu$ is approximated by hyperbolic sets in weak$*$-topology and in entropy (see for instance ~\cite[Proposition 1.4]{C}, ~\cite[Theorem 1]{Ge} and ~\cite{Ka}). Therefore, the main novelty here is to deal with non-hyperbolic ergodic measures whose existence is guaranteed by ~\cite{BBD2} (see also~\cite{BBD3,BZ1}). Before getting to that, let us recall the history of this topic. 

The topic of approximating invariant measures by   \emph{periodic measures} (i.e. atomic measures equidistributed on a single periodic orbit)
goes back to K. Sigmund~\cite{Sig1}: for hyperbolic basic sets, each invariant measure is approximated by  periodic measures  in weak$*$-topology. As it was shown in ~\cite{AS} that hyperbolic systems are not dense among differentiable diffeomorphisms. One can ask to what extent one can extend this result to non-hyperbolic setting. 

Under $C^{1+\alpha}$-setting,  many  results have already been obtained, for instance:
\begin{itemize} 
	\item A. Katok~\cite{Ka} has proven that each hyperbolic ergodic measure is approached by periodic measures in weak$*$-topology;
	\item In~\cite {SW}, the authors  have  shown that each  hyperbolic ergodic measure is  approximated by periodic measures in weak$*$-topology and in spectrum.
	\end{itemize} 

Under $C^1$-generic setting,   Ma\~n\'e's ergodic closing lemma shows that   every ergodic measure is approached by periodic measures in weak$*$-topology (see also ~\cite{ABC}). In  $C^1$-dominated setting, the result in ~\cite{C} states that hyperbolic ergodic measures are approximated by periodic measures in weak$*$-topology and the supports of hyperbolic ergodic measures are also approached by periodic orbits in Hausdorff topology at the same time (See also~\cite{G2} ). 

Beyond uniform hyperbolicity, the presence of non-hyperbolic ergodic measures is quite common, see~\cite{CCGYW}. However, there are few results on extending Sigmund's result to non-hyperbolic ergodic measures: ~\cite{DGR} shows that for some open set of step skew-products partially hyperbolic diffeomorphisms with compact  center foliation, each non-hyperbolic ergodic measure is approached by hyperbolic basic sets in weak$*$-topology and in entropy. The systems considered in ~\cite{DGR} exhibit the property called \emph{skeleton property}, and the skeleton property yields \emph{multi-variable-time-horseshoes} for each non-hyperbolic ergodic measure; finally such horseshoes give the hyperbolic sets approaching the non-hyperbolic ergodic measures in weah$*$-topology and in entropy.
In ~\cite{BZ}, for an open and dense subset of $\cU(M)$, the authors  show that each non-hyperbolic ergodic measure is approached by periodic measures of different indices in weak$*$-topology.  In~\cite{BZ}, the authors combine a shadowing lemma given by S. Liao~\cite{Liao1} with a tool called \emph{blender} designed by C. Bonatti and L. D\'iaz to get the sequence of periodic measures approaching non-hyperbolic ergodic measures. 

 K. Sigmund~\cite{Sig2} also initiated a topic to study the connectedness of the space of hyperbolic ergodic measures, and there was no further results after that. Until recently, A. Gorodetski and Y. Pesin~\cite{GP} show that  under $C^{1+\alpha}$-setting, the set of hyperbolic ergodic measures of a given index supported on an isolated homoclinic class (with some other properties) is \emph{path connected} and is \emph{entropy dense} among some certain set of hyperbolic invariant measures (see Theorems 1.1 and 1.5 in ~\cite{GP} for precisely statements).

One says that $\cM\subset \cM_{inv}(M,f)$ is \emph{path connected}, if for any $\mu,\nu\in\cM$, there exists a continuous path $\{\nu_t\}_{t\in[0,1]}$ in $\cM$ such that $\nu_0=\mu$ and $\nu_1=\nu$. 
Recall that given two subsets $\cM_1\subset \cM_2\subset \cM_{inv}(M,f)$, one says that $\cM_1$ is \emph{entropy dense} in $\cM_2$, if for any $\mu\in\cM_2$, there exists a sequence of measures $\{\mu_n \}\subset\cM_1$ such that $\mu_n$ tends to $\mu$ and $h_{\mu_n}(f)$ tends to $h_\mu(f)$.

As a corollary of Theorem~\ref{thmA}, one can also obtain the path connectedness and entropy density:
\begin{Co}\label{co:path-entropy}
For $f\in\cV(M)$ given by Theorem~\ref{thmA}, the set of ergodic measures is path connected and the set of hyperbolic ergodic measures is entropy dense in the set of ergodic measures.
\end{Co}
\begin{Remark} The path connected property is actually obtained in ~\cite{BZ} although it is not explicitly stated. 
	\end{Remark}
Combining Corollary~\ref{co:path-entropy} with the variational principle, one has that 
\begin{corollary} 
	For each $f\in\cV(M)$, one has that 
	$$h_{top}(f)=\sup\{h_\mu(f):\textrm{\:$\mu$ is a hyperbolic ergodic measure.}\}$$
	\end{corollary}
One can expect the following to be true:
\begin{Question*}
	Does there exist  a residual subset  of $\diff^r(M)$ $(r\geq 1)$ such that  the topological entropy equal the supremum of the metric entropy of hyperbolic ergodic measures?
\end{Question*}
Indeed, there are non-generic examples which  one can approach topological entropy by the metric entropy of non-hyperbolic ergodic measures, for instance the skew-products of the rotations on circles over the Anosov diffeomorphisms on rotus. As it is shown in ~\cite{BBD2,CCGYW}, there are rich families (open or generic) of  diffeomorphisms exhibitting non-hyperbolic ergodic measures; then  one can ask 
\begin{Question}
	Does there exist an open subset (or a locally generic subset) of $\diff^1(M)$ such that the topological entropy  equals  the supremum of the metric entropy of non-hyperbolic ergodic measures?
	\end{Question}
The results in ~\cite{TY} suggest that the answer might be negative (see also~\cite{DGR1}). In ~\cite{TY},  among certain $C^2$ open set of  partially hyperbolic (but non-hyperbolic) diffeomorphism with one-dimensional center, the authors show that  if the metric entropy of an ergodic measure is close enough to the topological entropy, then such ergodic measure is hyperbolic. 

{\bf Acknowlegement:}

We would like to thank Christian Bonatti and Sylvain Crovisier whose questions motivated this paper.


\section{Preliminary}
In this section, we collect the results and notions used in this paper.
\subsection{Entropy for horseshoes}
Let $\La$ be a hyperbolic horseshoe with $k$-legs, meaning that the dynamics in the horseshoe is conjugate to a full shift over $k$-symbols.  It is well known that $(\La,f)$ is conjugate to $(\Sigma_k,\sigma)$ and $h_{top}(f|_{\La})=\log k$.
\begin{theorem}\label{thm.middle-value-entropy}
	Let $\La$ be a hyperbolic horseshoe with $k$-legs. For any $r\in(0,\log k)$, there exists a compact invariant set $\La_r\subset\La$ such that $h_{top}(f|_{\La_r})=r$.
	\end{theorem}
The proof of  Theorem~\ref{thm.middle-value-entropy}  can be found in ~\cite[\S 7.3]{Wa}
\subsection{Partial hyperbolicity}
Let $K$ be a compact invariant set of a diffeomorphism $f$. Recall that $K$ is \emph{partially hyperbolic}, if the tangent bundle over $K$ is splitted into three invariant subbundles: one  is uniformly contracting, the other is uniformly expanding and the center has  intermediate behavior. To be precise, there exist an invariant splitting $T_KM=E^s\oplus E^c\oplus E^u$ and  a metric $\|\cdot\|$ such that for each $x\in K$, one has
$$\norm{Df|_{E^s(x)}}<\min\{1,m(Df|_{E^c(x)})\}\leq\max\{1,\norm{Df|_{E^c(x)}}\}<m(Df|_{E^u(x)}).$$

One says that $f$ is partially hyperbolic, if the whole manifold is a partially hyperbolic set. For a partially hyperbolic diffeomorphism, the strong stable and unstable bundles are uniquely integrable (see~\cite{HPS}), that is, there exist unique invariant  foliations $\cF^{ss}$ and $\cF^{uu}$ tangent to $E^s$ and $E^u$ respectively.  In this paper, the partially hyperbolic diffeomorphisms we consider have minimal strong stable and unstable foliations. We recall that a foliation is \emph{minimal}, if every leaf is dense in the manifold.
The following result is a folklore property for minimal foliations.
\begin{lemma}
Let $\cF$ be a minimal foliation on a compact manifold $M$. For any $\varepsilon>0$, there exists $L>0$ such that for any point $x$, one has that $\cF_{L}(x)$ is $\varepsilon$-dense in $M$, where $\cF_L(x)$ denotes the $L$-neighborhood of $x$ in the leaf $\cF(x)$.

\end{lemma}
It is interesting to notice the following corollary, although it is not used in this work.
\begin{corollary}
Let $f$ be a partially hyperbolic diffeomorphism. If one of the strong foliations is minimal, then  $f$ is topologically mixing.
\end{corollary}

Two hyperbolic periodic orbits are ~\emph{homoclinically related}, if the stable manifold of one intersects the unstable manifold of the other transversely, and vice versa. 
Here, we remark that for a partially hyperbolic diffeomorphism $f$ whose strong stable foliation is minimal, all the hyperbolic periodic orbits
of index $\dim(E^s)$ are homoclinically related.

\subsection{A shadowing lemma}

In this whole section, we assume that $\Lambda$ is a partially hyperbolic   set with one-dimensional center.
In this paper, we will use a   shadowing lemma  to find the periodic orbits near non-hyperbolic ergodic measures. The shadowing lemma  is firstly given by S. Liao~\cite{Liao1} and is improved by S. Gan~\cite{G1}, and  we will present a version adapted to  partially hyperbolic setting.
\begin{lemma}\label{l.shadowing}\cite{Liao1}
Let $\Lambda$ be a partially hyperbolic set with the splitting of the form $T_{\Lambda}M=E^s\oplus E^c\oplus E^u$  where $\dim(E^c)=1.$   For any $\lambda\in(0,1)$, there exist $L>0$ and $\e_0>0$ such that for any point $x\in\Lambda$,
 if there exists an integer $n$ such that
\begin{itemize}
\item[--] $$\ud(x,f^n(x))<\e\leq \e_0,$$
\item[--] $$\prod_{i=k}^{n-1}\norm{Df|_{E^c(f^{i}(x))}}\geq\lambda^{k-n}, \textrm{  for $k=0,\cdots,n-1$;}$$
\end{itemize}
then there exists a periodic point $p$ of period $n$ such that
$$\ud(f^i(p),f^i(x))<L\cdot \e, \textrm{ for $i=0,\cdots,n-1$}. $$
\end{lemma}

\subsection{Pliss lemma}
In this part,we recall the Pliss lemma ~\cite{P} which helps to find the hyperbolic time.

\begin{lemma}\label{l.pliss}
Given $a<b<c$,  there exists a constant $\rho>0$ such that for any family of real numbers $\{a_1,\cdots, a_n\}$ satisfying 
\begin{itemize}
	\item $a_i\leq c$ for each $i$;
	\item $\frac{1}{n}\sum_{i=1}^na_i\geq b$;
	\end{itemize}
 one has  $l$ integers $1\leq k_1<\cdots<k_l\leq n$ such that
\begin{itemize}
\item $l/n\geq \rho$;
\item for each $k_i$, one has  $\frac{1}{k_i-j+1}\sum_{s=j}^{k_i}a_s\geq a$, for any $j\in[1,k_i]$.
\end{itemize}
\end{lemma}
Recall that a continuous splitting $E\oplus F$ over the tangent bundle of a compact invariant set $\La$ is called \emph{dominated}, if there exists a metric $\norm{\cdot}$ such that 
$$\norm{Df|_{E(x)}}\cdot\norm{Df^{-1}|_{F(f(x))}}<\frac 12 \textrm{\:\: for any $x\in\La$}.$$
\begin{lemma}\label{l.unstable}
Let $\Lambda$ be a compact invariant set with dominated splitting $T_\La M=E\oplus F$. We fix a plaque family $W^{cu}$ associated to the bundle $F$. For any $\lambda\in(0,1)$, there exists $\delta>0$ such that for any point $x\in\Lambda$ satisfying  $$\prod_{j=0}^{n-1}\norm{Df^{-1}|_{F(f^{-j}(x))}}\leq \lambda^n\textrm{\:\: for any $n\in\mathbb{N}$},$$
one has that $W^{cu}_\delta(x)$ is contained in the unstable manifold of $x$.
\end{lemma}

\subsection{Blender}
\emph{Blender} is a powerful tool which is firstly  designed by C. Bonatti and L. D\'iaz to find robustly non-hyperbolic phenomenon by mixing dynamics of different behavior, see~\cite{BD1,BD2,BD3}. Recently, people start to use blender for finding robust phenomenon beyond non-uniform hyperbolicity, see for instance~\cite{BBD2, BBD3, BZ1}. There are several notions of a blender, see ~\cite{BD1,BD3, BBD2,BCDW}. In this section, we will recall the   main properties of a blender, and the blender  we will use in this paper  is a special blender called \emph{blender horseshoe} ~\cite{BD3}. 

A $cu$-blender horseshoe $\La$ is a hyperbolic basic set and $T_\La M$ is splitted into dominated subbundles of the form $T_{\La}M=E^s\oplus E^{cu}\oplus E^{u}$ such that 
\begin{itemize}
	\item $E^{cu}$ is one dimensional and is uniformly expanding;
	\item $\La$ is the maximal invariant set in an embedded cube $K=[-1,1]^s\times[-1,1]\times[-1,1]^u$, where $s=\dim(E^s), u=\dim(E^u)$;
	\item the dynamics in $K$ is a horseshoe with 2-legs;
	\item there exist cone fields $\cC^s,\cC^u,\cC^{uu}$ defined in $K$ such that the tangent space of $[-1,1]^s\times\{y\}\times\{z\}$ is contained in $\cC^s$, the tangent space of $\{x\}\times[-1,1]\times[-1,1]^u$ is contained in $\cC^u$ and the tangent space of $\{x\}\times\{y\}\times[-1,1]^u$ is contained in $\cC^u$.
	\item $\cC^s$ is $Df^{-1}$-invariant, and $\cC^u,\cC^{uu}$ are $Df$-invariant;
	\item $Df$ is uniformly expanding along the cone field $\cC^u$, that is, there exists $\tau>1$ such that for any $x\in K\cap f^{-1}(K)$ and any unit vector $v\in\cC^{u}(x)$, one has $\norm{Df(v)}\geq \tau.$
	\item A disc $D$ tangent to the cone field $\cC^{uu}$ is called $uu$-disc, if the boundary of $D$ is contained in $[-1,1]^s\times[-1,1]\times\partial{([-1,1]^u)}$ and the interior of $D$ is contained in the interior of $K$; there exists an open region in $K$ called the \emph{superposition region} of $\La$ such that for any $uu$-disc $D$ intersecting this region, at least one of the connected components of  $f(D)\cap K$  is a $uu$-disc intersecting the open region. 
\end{itemize}
See ~\cite[Section 3.2]{BD3} for the precise definition of a blender horseshoe. We will call the cube $K$ as the \emph{characteristic region} of the blender horseshoe $\La$. We remark that the existence of a blender horseshoe is an open property, that is, there exists a small neighborhood $\cU$ of $f$ such that for $g\in\cU$, the cotinuation $\La_g$ of $\La$ is a blender horseshoe.

A disc $S$ of dimension $\dim(E^{cu}\oplus E^{u})$ is called a $cu$-strip, if there exists a $C^1$-embedding map $\phi:[-1,1]\times[-1,1]^u\mapsto K$ such that 
\begin{itemize}
	\item $\phi([-1,1]\times[-1,1]^u)=S$
	\item  for any $t\in[-1,1]$, the disc $\phi(\{t\}\times [-1,1]^u)$ is a $uu$-disk;
	\item $S$ is disjoint from $\partial([-1,1]^s)\times[-1,1]\times[-1,1]^u.$
\end{itemize}
The  $uu$-discs $\phi(\{-1\}\times[-1,1]^u)$ and $\phi(\{1\}\times[-1,1]^u)$ are called the \emph{vertical boundary of $S$}. The \emph{center length $\ell^c(S)$ of $S$} is defined as the minimum length among all the $C^1$ curves in $S$ which joins  the two vertical boundary of $S$. Hence, each $cu$-strip $S$ contains a disc tangent to $\cC^u$ of radius at least $\ell^c(S)/2$.

One important property of a $cu$-strip is the following result.
\begin{lemma}\label{l.iterate-cu-strip}\cite[Lemma 3.20]{BZ}
	There exists $c>0$ such that for any $cu$-strip  $S$ satisfying 
	\begin{itemize}
		\item $\ell^c(S)<c$;
		\item $S$ contains a $uu$-disc intersecting the superposition region,
	\end{itemize}
	the disc $f(S)$ contains a $cu$-strip $S_1$ such that 
	\begin{itemize}
		\item[--] $\ell^c(S_1)>\tau\cdot\ell^c(S)$;
		\item[--] $S_1$ contains a $uu$-disc intersecting the superposition region.
	\end{itemize}
\end{lemma} 

One says that a diffeomorphism  exhibits a \emph{co-index one heterodimensional cycle}, if there exist hyperbolic periodic orbits $P$ and $Q$ such that $Ind(P)=Ind(Q)+1$ such that $W^s(P)$ has non-empty transverse intersection with $W^u(Q)$ and $W^u(P)$ intersects $W^s(Q)$.
In~\cite{BD3}, the authors show that a co-index one heterdimensional cycle can yield a $cu$-blender horseshoe. 
\begin{theorem}\label{thm.existence-pf-blender}~\cite[Proposition 5.6]{BD3}
	Let $f$ be a $C^1$ diffeomorphism with a co-index one  heterdimensional cycle. Then there exists $g$ arbitrarily $C^1$ close to $f$ such that $g$ has a $cu$-blender horseshoe $\La_g$.
\end{theorem}

\subsection{Metric entropy defined by A. Katok}\label{s.entropy-katok}
In this part, we recall the definition of metric entropy of an ergodic measure given by A. Katok~\cite{Ka}.

Let $f:X\mapsto X$ be a homeomorphism on the compact metric space $(X,\ud)$. For $\varepsilon>0$ and $n\in\N$, a \emph{$(n,\varepsilon)$-Bowen ball centered at $x$} is the open set $$B(x,\varepsilon)\cap f^{-1}(B(f(x),\varepsilon))\cap\cdots\cap f^{-n}(B(f^n(x),\varepsilon)),$$
where $B(y,\e)$ denotes the $\e$-ball centered at $y\in X$.

Let $\mu$ be an ergodic measure of $f$. For $\delta>0$ and $\varepsilon>0$, let $N(n,\varepsilon,\delta)$ be the minimal  number of $(n,\varepsilon)$-Bowen balls which cover a set of measure no less than $1-\delta$.  Then one has 
\begin{theorem}\label{thm.entropy-katok}~\cite[Theorem I.I]{Ka} 
	Let $\mu$ be an ergodic measure, then for $\delta\in(0,1)$  one has
$$h_{\mu}(f)=\lim_{\varepsilon\rightarrow 0}\liminf_{n\rightarrow\infty}\frac{\log{N(n,\varepsilon,\delta)}}{n}=\lim_{\varepsilon\rightarrow 0}\limsup_{n\rightarrow\infty}\frac{\log{N(n,\varepsilon,\delta)}}{n}.$$
\end{theorem}

Recall that 
a set $K_1$ is called a \emph{$(n,\e)$-separated set} of $K_2$, if $K_1\subset K_2$ and for any two different points $x,y\in K_1$, there exists $j\in\{0,\cdots, n\}$ such that $\ud(f^j(x),f^j(y))\geq\e$.
As a corollary of Theorem~\ref{thm.entropy-katok}, one has that 
\begin{corollary}\label{co:entropy-separated-set}
	Let $\mu$ be an ergodic measure and $K$ be a measurable set with $\mu(K)\geq 1-\delta$   for some $\delta\in(0,1)$. We denote by $S(n,\varepsilon,K)$ be the maximal number of $(n,\varepsilon)$-separated sets of $K$, then one has 
	$$\lim_{\varepsilon\rightarrow 0}\limsup_{n\rightarrow\infty}\frac{\log{S(n,\varepsilon,K)}}{n}\geq h_\mu(f).$$
	\end{corollary}

\section{Periodic orbits near non-hyperbolic ergodic measures}

One of the  ingredients for the proof of Theorem~\ref{thmA} is the following result, and some ideas can be found in  ~\cite{BZ}.

Recall that each diffeomorphism  $f\in\cU(M)$  is partially hyperbolic with one dimensional center and exhibits periodic points of different indices; moreover,  the strong foliations of $f$ are robustly minimal. The minimality of strong foliations is fully used in this section.

\begin{theorem}\label{thm.shadow}
There exists a $C^1$ open and dense subset $\cV(M)$ of $\cU(M)$ such that for  $f\in\cV(M)$,
there exists a  constant $c_0>0$  such  that for any $\varepsilon>0$, $\chi>0$  small and  any $C>1$,  there exist two positive  integers $N(\varepsilon)\ll N(C,\chi,\varepsilon)=N $ such that
for any $x\in M$ satisfying 
$$\frac{1}{C}\cdot e^{-k\chi}\leq\norm{Df^k|_{E^c(x)}}\leq C\cdot e^{k\chi}\textrm{\:\:for any $k\in\mathbb{N}$},$$
and  any $n>N$, there exist $z\in M$ and $m\in(n,n+c_0\cdot n\cdot \chi)$ such that
\begin{itemize}
\item[--] $$\ud(z,f^m(z))<\varepsilon$$
\item[--]  $$\ud(f^i(f^{N(\varepsilon)}(z)),f^i(x))<\varepsilon \textrm{ for $i=0,\cdots,n-N(\e)$;}$$
\item[--]$$\frac{1}{m-j}\sum_{i=j}^{m-1}\log{\norm{Df|_{E^c(f^i(z))}}}\geq  4\chi \textrm{\: for $j=0,\cdots,m-1$}.$$
\end{itemize}
\end{theorem}
 
Before giving the proof, let us explain the idea. Firstly, we show that for an open and dense subset in $\cU(M)$, there exist $cu$-blender horseshoes. One fixes a hyperbolic periodic point $p$ with positive center Lyapunov exponent and  a point $x$ satisfying the properties in the assumption of Theorem~\ref{thm.shadow}, then the minimality of strong stable foliation ensures that there exists $y\in W^u(p)\cap \cF^{ss}(f^{-N}(x))$ (for some $N$) that is arbitrarily close to $p$.  We will choose a small disc  $D$ tangent to $E^c\oplus E^u$ such that it is contained in local unstable manifold of $p$ and contains $y$. We will iterate  $D$ forwardly to follow the positive orbit of $y$ for a long time (say time $n$) and therefore can follow the orbit of $x$ for a long time (almost time $n$), then by the minimality of strong unstable foliation, after a small proportion  of time (with the respect to $n$), the images of $D$ contain a $cu$-strip in the superposition region of the $cu$-blender. The property of a cu-blender horseshoe guarantees that one can iterate the $cu$-strip for a proportion of time (about $n\chi$-level) to ensure that the $cu$-strip has some uniform size. Once again, by the minimality of the strong stable foliation, the fixed size of  strong stable manifold of $p$ intersects the $cu$-strip with uniorm size,  and this can yield the orbit segment with posited properties. 
\proof[Proof of Theorem~\ref{thm.shadow}]
By the robust minimality of strong foliations, the diffeomorphisms in $\cU(M)$ are robustly transitive. Since there are periodic points of different indices, by  Hayashi's connecting lemma~\cite{H}, there exist a $C^1$ dense subset of $\cU(M)$ such that each diffeomorphism in this dense subset exhibits a co-index one heterodimensional cycle (for details see for instance~\cite[Lemma 2.5]{BDPR}). By Theorem~\ref{thm.existence-pf-blender}, there exists a $C^1$ dense subset of $\cU(M)$ exhibitting $cu$-blenders  whose existence is a robust property. To summarize,   there exists a $C^1$ open and dense subset $\cV(M)$ of $\cU(M)$ such that each $f\in \cV(M)$ has a $cu$-blender horseshoe $\Lambda^u$.

From now on, we fix $f\in\cV(M)$, a hyperbolic periodic point $p$ of $f$ whose center Lyapunov exponent is positive and a $cu$-blender horseshoe $\La^u$ of $f$.  Without loss of generality, we assume that $p$ is a fixed point. We denote by $2\lambda$ the center Lyapunov exponent of $p$, then we fix its  local stable and unstable manifolds  $W^s_{2\varepsilon_0}(p)$ and $W^u_{2\varepsilon_0}(p)$   for some $\varepsilon_0>0$ small  such that for any $w\in B_{2\varepsilon_0}(p)$, one has
$$\log\norm{Df^{-1}|_{E^c(w)}}\leq -\frac{3\lambda}{2}.$$
Denote by $K$ the characteristic region of the blender $\La^u$. By definition, there exists $\tau>1$ such that 
 $$ \norm{Df|_{E^c(w)}}\geq \tau \textrm{\:\: for $w\in K\cap f^{-1}(K)$}.$$
 The posited constant $c_0>0$ would be given by $\lambda$ and $\tau$, and the precise formula would be given in the end.
 
 \medskip
 
For $\varepsilon,\chi\in(0,\min\{\frac{\lambda}{100},\frac{\log\tau}{100},\e_0\})$ small and $C>1$, we will find the posited integers $N(\e)$ and $N(\chi,C,\e)$.  

  By the minimality of the strong stable foliation $\cF^{ss}$, there exists    $l=l(\varepsilon)>0$ such that for any point $w\in M$, the strong stable disc $\cF^{ss}_{l}(w)$ intersects  $W^u_{\varepsilon/4}(p)$.
  By the uniform contraction of $f$ along $\cF^{ss}$,
    there exists an   integer $N(\varepsilon)\in\mathbb{N}$ such that $$f^{N(\varepsilon)}(\cF^{ss}_{l}(w))\subset \cF^{ss}_{\varepsilon/2}(f^{ N(\varepsilon)}(w))\textrm{\: for any $w\in M.$}$$
Let $x\in M$ be a point satisfying 
$$\frac{1}{C}\cdot e^{-k\chi}\leq\norm{Df^k|_{E^c(x)}}\leq C\cdot e^{k\chi}\textrm{\:\:for any $k\in\mathbb{N}$}.$$ By the choice of $l$, there exists a transverse intersection $y\in  \cF^{ss}_{l}(f^{-N(\varepsilon)}(x))\cap W^u_{\varepsilon/4}(p)$. By the choice of $y$,  one has 
\begin{itemize} 
	\item  $f^{N(\e)}(y)\in \cF^{ss}_{\e/2}(x)$;
	\item $y$ has unstable manifold of size $\varepsilon_0$ contained in $W^u_{2\varepsilon_0}(p)$.
	\end{itemize}By the uniform contraction of $f$ along $\cF^{ss}$ and the uniform continuity of the center bundle $E^c$, 
 there exists $C^{\prime}=C^{\prime}(C,\chi,\e)>1$ which is independent of $x$ such that  
$$\frac{1}{C^{\prime}}\cdot e^{-2k\chi}\leq\norm{Df^n|_{E^c(y)}}\leq C^{\prime}\cdot e^{2k\chi},\textrm{ for any $k\in\mathbb{N}$}.$$

 To continue the proof, we need a new notation. 
We fix a $Df$-strictly invariant cone field $\cC^{uu}$ of $E^{u}$ along which $Df$ is uniformly expanding. Since $E^c$ is transverse to $E^u$ restricted to $E^c\oplus E^u$ and the angle is uniformly bounded from below, for simplicity, one can assume that $E^c$ is orthorgonal to $E^u$.  We denote by $u=\dim(E^{u})$. 
One says that a disc $D$ tangent to $E^c\oplus E^u$ is \emph{$C^1$-foliated by discs tangent to $\cC^{uu}$}, if there exists a $C^1$ embedding
$\varphi:[0,1]\times [0,1]^u\mapsto M$ such that
\begin{itemize}
	\item[--] $\varphi([0,1]\times [0,1]^u)=D$;
	\item[--] $\varphi(\{t\}\times[0,1]^u)$ is a disc tangent to $\cC^{uu}$ for $t\in[0,1]$;
\end{itemize}
The discs $\varphi(\{0\}\times[0,1]^u)$ and $\varphi(\{1\}\times[0,1]^u)$ are called the \emph{vertical boundary} of the $cu$-disc $D$.
If there exists a $C^1$ curve tangent to $E^c$  contained in $D$ whose two endpoints are  in the interior of $\varphi(\{0\}\times[0,1]^u)$ and $\varphi(\{1\}\times[0,1]^u)$ respectively, then we will  call $D$ a \emph{vertical $cu$-disc}.
For each vertical  $cu$-disc $D$, the \emph{center length} $\ell^c(D)$ of $D$ is defined as the infimum length of center curves in $D$ with two endpoints contained in the interior of $\varphi(\{0\}\times[0,1]^u)$ and $\varphi(\{1\}\times[0,1]^u)$ respectively.

By the uniform continuity of $E^c$,  there exists $\delta\in(0,\e/4)$ such that
for any $w_1,w_2\in M$ with $\ud(w_1,w_2)<\delta$, one has that
$$|\log\norm{Df|_{E^c(w_1)}}-\log\norm{Df|_{E^c(w_2)}}|<\chi.$$

 For $n\in\mathbb{N}$, we consider the $\delta\cdot e^{-5n\chi}$-tubular neighborhood of $W^{uu}_{\delta}(y)$ in $W^{u}_{\varepsilon_0}(y)$ and we denote it as $S_n$. Then by the transversality between $E^c$ and $E^u$,  there exists an integer  $N_0(\chi)$ which depends only on $\chi$ such that for $n\geq N_0(\chi)$, one has 
 \begin{itemize}
 	\item $S_n$ is foliated by discs tangent to $\cC^{uu}$;
 	\item the length of each center curve through $y$ contained in $S_n$ is bounded by $\delta\cdot e^{-4n\chi}$;
 	\item each center curve through $y$ and contained in $S_n$ would intersect the two  vertical boundaries of $S_n$  into unique points whose distance to $y$ is less than $\delta/4$.
 	\end{itemize}
 Now, we will iterate $S_n$ forwardly, and each time we cut it by the $\delta$-neighborhood of the forward orbit of $y$. To be precise,  for any $i\in\mathbb{N}$, we consider the connected component $S_n(i)$ of $f^i(S_n)\cap B_{\delta}(f^i(y))$ containing $f^i(y)$.  By the strictly invariant property of the cone field $\cC^{uu}$, $S_n(i)$ is $C^1$ foliated by discs tangent to $\cC^{uu}$.
\begin{lemma}\label{l.iterate} There exists an integer $N_1=N_1(C, \chi,\e)$ such that for  $n>N_1$, one has that for $i\leq n$, the cu-disc $S_n(i)$ is a vertical $cu$-disc satisfying that 
	\begin{itemize}
		\item $\delta\cdot e^{-9n\chi}\leq \ell^c(S_n(i));$
		\item any center curve contained in $S_n(i)$ and through $f^i(y)$ intersects the two vertical boundaries of $S_n(i)$ into unique points whose distance to $f^i(y)$ is less than $\frac{\delta}{2}$. 
		\end{itemize}
	\end{lemma}
\proof 

For any $C^1$-curve $\gamma$, we denote by $\ell(\gamma)$ the length of   $\gamma$.

\begin{claim}\label{c.first center} 
	There exists an integer $N_1=N_1(C,\chi,\e)>N_0$ such that for $n\geq N_1$ and any center curve $\gamma$ through $y$  whose endpoints are contained in the vertical boundary of $S_n$, one has that
	$$\delta\cdot e^{-9n\chi}\leq \ell(f^i(\gamma))<\frac{\delta}{2} \textrm{\: for any $i\in\{0,1,\cdots,n \} $}.$$
\end{claim}
\proof
Let $N_1=N_1(C,\chi,\e)$ be the smallest integer such that 
  $$N_0<N_1\textrm{\: and\:} C^{\prime}\cdot e^{-N_1\chi}<\frac{1}{2}.$$
 
Let $\gamma$ be a center curve through $y$   with endpoints contained in the vertical boundary of $S_n$. Then $\ell(\gamma)\leq \delta\cdot e^{-4n\chi}<\delta/2.$
If $\ell(f^j(\gamma))<\delta/2$ for $j\leq i$, then by the choice of $\delta,$ one has the estimate:
\begin{align*}
\ell(f^{i+1}(\gamma))&=\int\norm{Df^{i+1}(\gamma^{\prime}(t))}\ud t\\
&=\int\norm{Df^{i+1}|_{E^c(\gamma(t))}}\norm{\gamma^{\prime}(t)}\ud t\\
&\leq  e^{(i+1)\chi} \cdot \norm{Df^{i+1}|_{E^c(y)}}\cdot \ell(\gamma)\\
&\leq e^{(i+1)\chi} \cdot
C^{\prime} \cdot e^{2(i+1)\chi}\cdot \delta\cdot e^{-4n\chi}\\
&\leq\frac{\delta}{2}.
\end{align*}
Now we can give the lower bound of  $\ell(f^i(\gamma))$ for $i\leq n$. Since $\ell(f^i(\gamma))<\delta/2$, one has 
\begin{align*}
\ell(f^{i}(\gamma))&=\int\norm{Df^{i}|_{E^c(\gamma(t))}\gamma^{\prime}(t)}\ud t\\
&=\int\norm{Df^{i+1}|_{E^c(\gamma(t))}}\norm{\gamma^{\prime}(t)}\ud t\\
&\geq  e^{-(i+1)\chi} \cdot \norm{Df^{i+1}|_{E^c(y)}}\cdot \ell(\gamma)\\
&\geq  e^{-(i+1)\chi}\cdot \frac{1}{C^{\prime}} \cdot e^{-2(i+1)\chi}\cdot \delta\cdot e^{-5n\chi}\\
& \geq \delta\cdot e^{-9n\chi}.
\end{align*}

\endproof
 
By Claim~\ref{c.first center}, for $n\geq N_1$, each $S_n(i)$ contains center curves through $y$ and whose endpoints are contained in the vertical boundary of $S_n(i)$, then by definition, the $cu$-disc $S_n(i)$ is a vertical $cu$-disc. Morover, any center curve contained in $S_n(i)$ and through $f^i(y)$ intersects the two vertical boundaries of $S_n(i)$ into unique points whose distance to $f^i(y)$ is less than $\delta/2$.

Now, we give a lower bound for $\ell^c(S_n(i))$. 
\begin{claim} For $i\leq n$, the center length of $S_n(i)$ is lower bounded by $\delta\cdot e^{-9n\chi}$, in formula:
$$\ell^c(S_n(i))\geq \delta\cdot e^{-9n\chi}.$$
\end{claim}
\proof
By definition, the center length  of  $S_n(0)=S_n$ is lower bounded by $\delta\cdot e^{-5n\chi}$ and the center length  of $S_n(0)$ is upper bounded by $\delta$. We will prove the claim by induction. Assume that the claim is true  for $j\leq i$.  For any center curve $\gamma\subset S_n(i+1)$ whose endpoints are contained in the vertical boundary of $S_n(i+1)$, by the uniform expansion of $Df$ along $\cC^{uu}$,
one has that
\begin{itemize}
\item[--] $f^{-j}(\gamma)$ is contained in $S_n(i+1-j)$,  whose endpoints are contained in the vertical boundary of $S_n(i+1-j)$, for $j\leq i+1$;
\item[--] the center segment $f^{-j}(\gamma)$ is contained in the $\delta$-neighborhood of $f^{i+1-j}(y)$.
    \end{itemize}
By the choice of $\delta$, one has the estimate:
\begin{align*}
\ell(\gamma)&=\int\norm{Df^{i+1}|_{E^c(f^{-i-1}(\gamma(t)))}}\cdot\norm{\frac{\ud f^{-i-1}(\gamma(t))}{\ud t}}\ud t\\
&\geq \frac{1}{C^\prime}\cdot e^{-2(i+1)\chi} \cdot e^{-(i+1)\chi} \cdot \ell(f^{-i-1}(\gamma))\\
&\geq \delta\cdot  e^{-9n\chi}.
\end{align*}
\endproof

Now, the proof of Lemma~\ref{l.iterate} is complete. \endproof

By the minimality of the strong unstable foliation $\cF^{uu}$, there exists $N_2=N_2(\delta)>0$ such that for any strong unstable disc $D$ of diameter at least $\delta/2$, one has that $f^{N_2}(D)$ contains a strong unstable disc in the superposition region of the blender horseshoe $\La^u$. By definition of $S_n(n)$, one has that  $\cF^{uu}_\delta(f^n(y))$ is contained in $S_n(n)$. By the uniform expansion of $Df$ along the cone field $\cC^{uu}$ and the second item in Lemma~\ref{l.iterate}, one has that  $f^{N_2}(S_n(n))$ contains a $cu$-strip $S$ whose center length is at least $\delta\cdot e^{-9n\chi}\cdot m_1^{N_2}$, where $m_1=\inf_{w\in M}\norm{Df|_{E^c(w)}}$. Now, the property of the $cu$-blender horseshoe $\La^u$ given by Lemma~\ref{l.iterate-cu-strip} allows us to iterate $S$ in the characteristic region of $\La^u$ until we get a $cu$-strip with uniform center length. By Lemma~\ref{l.iterate-cu-strip}, there exist a  constant $\e_1>0$ which depends only on $\La^u$ and an integer $N_3\in\mathbb{N}$ such that
\begin{itemize}
\item[--]$$N_3\leq \frac{9n\chi-N_2\log{m_1}+\log\delta+\log\e_1}{\log\tau};$$
\item[--] $f^{N_3}(S)$ contains a disc $D_n$ tangent to $E^{cu}$ of diameter at least $\varepsilon_1$.
\end{itemize}
Now, by the minimality of the strong stable foliation, there exists an integer $N_4=N_4(\varepsilon_1)$ such that $f^{-N_4}(W^{s}_{\varepsilon_0}(p))$ intersects $D_n$ in the interior (remember that $W^s(p)=\cF^{ss}(p))$), and we denote $\tilde{z}$ as the intersection. Let $z=f^{-n-N_2-N_3}(\tilde{z})$, then one has

\begin{itemize}
	\item[--] $z\in W^u_{\varepsilon/4}(p)$ and $f^{n+N_2+N_3+N_4}(z)\in W^s_{\e_0}(p)$;
\item[--] $\ud(f^i(z),f^i(y))<\delta<\frac{\e}{4}$, for $i\in[0,n]$;
\item[--] $f^j(z)$ is in the characteristic  region of $\La^u$, for $j\in[n+N_2, n+N_2+N_3]$.
\end{itemize}
Since $f^{N(\e)}(y)\in\cF^{ss}_{\e/2}(x)$, one has that   $\ud(f^i(f^{N(\varepsilon)}(z)),f^i(x))<\varepsilon$ for $i\in[0,n-N(\varepsilon)]$.
By the choice of $\delta$, one has
 $$\log{C^\prime}+3n\chi\geq\sum_{i=0}^{n-1}\norm{Df|_{E^c(f^i(z))}}\geq -\log{C^\prime}-3n\chi;$$

Now, we choose $t=[\frac{10\chi}{\lambda}\cdot n]+1$.
Let  $m=n+t+N_2+N_3+ N_4$.

\begin{claim}There exists an integer  $N=N(C,\varepsilon,\chi)>N_1+N(\e)$   such that for $n\geq N$, one has
\begin{itemize}
\item $f^m(z)\subset W^s_{\varepsilon/2}(p)$.
\item $$\frac{t+N_2+N_3+ N_4}{n}< \big(\frac{10}{\log\tau}+\frac{10}{\lambda}\big)\cdot\chi,$$

\item $$\frac{1}{m-j}\sum_{i=j}^{m-1}\log{\norm{Df|_{E^c(f^i(z))}}}\geq  4\chi \textrm{\: for $j=0,\cdots,m-1$}.$$

\end{itemize}
\end{claim}

 \begin{figure}
	\begin{center}
		\def\svgwidth{1\columnwidth}
		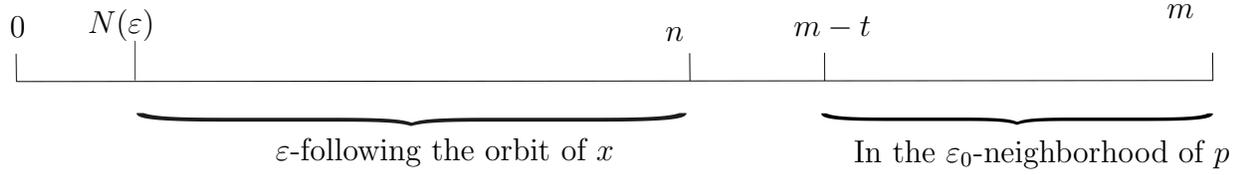
		\caption{The itinerary of $z$}
	\end{center}
\end{figure}

\proof

Since $t= [\frac{10\chi}{\lambda}\cdot n]+1$, $W^s(p)=\cF^s(p)$ and the choice of $N_2+N_4+N(\e)$ is independent of $n$,  therefore for $n$ large, the first and the second items are satisfied. 

Recall that the center Lyapunov exponent of $p$ is $2\lambda$ and $f^{n+N_2+N_3+N_4}(z)\in W^s_{\e_0}(p)$, then by the choice of $\varepsilon_0$, one has that
$$\sum_{i=m-j}^{m-1}\log\norm{Df|_{E^c(f^i(z))}}\geq j\cdot\frac{3\lambda}{2},\textrm{\: for $j\leq t$}.$$
Once again, since the choice of $N_4$ is independent of $n$,  for $n$ large, one has
 $$\frac{1}{j}\sum_{i=m-j}^{m-1}\log\norm{Df|_{E^c(f^i(z))}}\geq \lambda\textrm{\:\:\: for $j\in[t,t+N_4]$};$$

Now, we estimate the average along  whole orbit segment.

\begin{align*}
\sum_{i=0}^{m-1}\log{\norm{Df|_{E^c(f^i(z))}}}&=\sum_{i=0}^{n-1}\log{\norm{Df|_{E^c(f^i(z))}}}+\sum_{i=n}^{n+N_2-1}\log{\norm{Df|_{E^c(f^i(z))}}}
\\
&\hspace{5mm}+\sum_{i=n+N_2}^{n+N_2+N_3-1}\log{\norm{Df|_{E^c(f^i(z))}}}+\sum_{i=n+N_2+N_3}^{m-1}\log{\norm{Df|_{E^c(f^i(z))}}}
\\
&\geq -\log{C^{\prime}}-3n\chi+ N_2\cdot \log m_1+ N_3\cdot\log\tau +(t+N_4)\lambda.
\end{align*}
Since $\chi$ has been chosen in $(0,\min\{\frac{\lambda}{100},\frac{\log\tau}{100},\e_0\})$, for $n$ large, one has 
\begin{align*}
\sum_{i=0}^{m-1}\log{\norm{Df|_{E^c(f^i(z))}}}&\geq 
 -4n\chi  +(t+N_4)\lambda
\\
&=m\cdot\bigg(-4\chi\cdot\frac{n}{n+t+N_2+N_3+N_4} +\frac{t}{n+t+N_2+N_3+N_4}\lambda\bigg)
\\
&\geq m\cdot \frac{\frac{t}{n}\lambda-4\chi}{1+\frac{t}{n}+\frac{N_2+N_3+N_4}{n}} 
\\
&\geq m\cdot \frac{6\chi}{1+\frac{10\chi}{\lambda}+\frac{10\chi}{\log\tau}}
\\
&\geq m\cdot 5\chi.
\end{align*}

Let $N_5=N_5(C,\chi,\e)$ be an integer such that for any $j\geq N_5$, one has 
$$\log C^\prime+3\cdot j\cdot\chi<\frac{7}{2}\cdot j\cdot\chi.$$

Let $\rho$ be the number given by Lemma~\ref{l.pliss} corresponding to $4\chi<5\chi<\log m_2$, where $m_2=\sup_{w\in M}\norm{Df|_{E^c(w)}}$.

For $n$ large, one has that $$\frac{N_1+N_2+N_4+N_5}{n}\leq \frac{\rho}{2}.$$

Recall that for any $w\in K\cap f^{-1}(K)$ in the characteristic region of the blender $\La^u$. one has
$$\log \norm{Df|_{E^c(w)}}\geq \log\tau>5\chi.$$

By Lemma~\ref{l.pliss} and the fact that  $$\frac{1}{j}\sum_{i=m-j}^{m-1}\log\norm{Df|_{E^c(f^i(z))}}\geq \lambda>5\chi\textrm{\:\:\: for $j\in[t,t+N_4]$},$$
one has $$\frac{1}{m-j}\sum_{i=j}^{m-1}\log{\norm{Df|_{E^c(f^i(z))}}}\geq  4\chi \textrm{\: for $j=0,\cdots,m-1$}.$$
\endproof

One takes  $c_0=\frac{10}{\log\tau}+\frac{10}{\lambda}$, and this ends the proof of Theorem~\ref{thm.shadow}. 
\endproof

\section{A sufficient condition for horseshoes with large entropy}
We divide this section into two subsections. In the first part, we   show that in the dominated setting,  certain separated sets in a Pesin block  can give horseshoes  with large entropy.
 In the second part, we   show that for each non-hyperbolic ergodic measure, one can find Pesin block `close to' it;  moreover in the Pesin block, one can find some separated sets that can be used to build horseshoes.
 
 Given a compact manifold $M$, we denote by $\cM$ the set of probability measures on $M$. We fix a metric on $\cM$ defined by
 $$\ud(\mu,\nu)=\sum_{n\in\mathbb{N}}\frac{|\int\varphi_n\ud\mu-\int\varphi_n\ud\nu|}{2^{n+1}\cdot\|\varphi_n\|}\textrm{\: for $\mu,\nu\in\cM$,}$$
 where $\{\varphi_n\}_{n\in\mathbb{N}}$ is a dense subset of the set of continuous functions on $M$. 
 By the choice of the metric, for $\mu,\nu\in\cM$, one has $\ud(\mu,\nu)\leq 2.$

\subsection{Horseshoes derived from separated sets}
For a periodic point $p$, we denote by $\cO_p$ the orbit of  $p$ and by $\delta_{\cO_p}$ the periodic measure on $\cO_p$. Let $\{p_1,\cdots,p_l\}$ be a set of periodic points of $f$, then one can define a subset of invariant measures as follows:
$$\Cov(p_1,\cdots,p_l):=\big\{\sum_{i=1}^l\alpha_i\delta_{\cO_{p_i}}: \alpha_i\geq 0\: \textrm{ and }\: \sum_{i=1}^l\alpha_i=1 \big \}.$$

The aim of this section is to prove the following theorem.

\begin{theorem}\label{thm.entropy-periodic-orbit}
Let $f\in\diff^1(M)$ and $\La$ be a compact invariant set with the dominated splitting of the form $T_\La M=E\oplus F$. For  $\lambda\in(0,1)$, consider the Pesin block of the following form
$$\La_\lambda=\big\{x\in\La: \prod_{i=0}^{n-1}\norm{Df|_{E(f^i(x))}}\leq \lambda^n,\:\:  \prod_{i=0}^{n-1}\norm{Df^{-1}|_{F(f^{-i}(x))}}\leq \lambda^n, \textrm{ for $n\in\N$}\big\}.$$
For any $\varepsilon>0$, there exists $\eta_0>0$ such that for any $n\in\mathbb{N}$,   any sequence of periodic points $\{p_1,\cdots,p_k\}\subset \La_\lambda$ of  period  $n$ and any $\eta<\eta_0$,
if one has that
\begin{itemize}
\item[--] $\ud(p_i,p_j)<\eta$, for $1\leq i<j\leq k$;
\item[--] $\{p_1,\cdots,p_k\}$ is a $(n,16\cdot\eta)$ separated set.
\end{itemize}
Then there exists a hyperbolic set $K_{\varepsilon}$ of index $\dim(E)$ such that 
\begin{itemize}
	\item[--]
$$h_{top}(f|_{K_{\varepsilon}})\geq \frac{\log k }{n}.$$
\item[--] for any   $\nu\in\cM_{inv}(K_{\e},f)$, the measure $\nu$ is in the $\e$-neighborhood of $\Cov(p_1,\cdots,p_k)$.
\end{itemize}
\end{theorem}
\proof  Since $E$ is transverse to $F$, for simplicity, we will assume that $E$ is orthogonal to $F$. We fix plaque families $W^{cs}$ and $W^{cu}$ corresponding to $E$ and $F$ respectively. By Lemma~\ref{l.unstable},  there exists $\eta_\lambda>0$ such that for each $x\in\La_\lambda$, $W^{cs}_{\eta_\lambda}(x)$ and $W^{cu}_{\eta_\lambda}(x)$ are contained in the stable and unstable manifolds of $x$ respectively.

For $\varepsilon\in(0,\min\{\frac{\lambda}{4},\frac{1-\lambda}{4}\})$ small, there exist $\eta_0\in(0,\frac{\eta_\lambda}{4})$ and  two cone fields $\cC^E$ and $\cC^F$ defined in the $8\eta_0$-neighborhood of $\La$ such that
\begin{itemize}
		\item for any $x,y$ with $\ud(x,y)\leq 8\eta_0$, one has $\ud(\delta_x,\delta_y)<\e.$
		
	\item $\cC^E$ is $Df^{-1}$-invariant and $\cC^F$ is $Df$-invariant;
	
	\item  for any $x,y\in B_{8\eta_0}(\La)$ with $\ud(x,y)<8\eta_0$, any discs $D_x^E$ and $D^E_y$ through the points $x, y$  and  tangent to $\cC^E$ respectively, one has that $Df|_{T_xD_x}$ is $\varepsilon$ close to $Df|_{T_yD_y}$ and $Df^{-1}|_{T_xD_x}$ is $\varepsilon$ close to $Df^{-1}|_{T_yD_y}$ ; the analogous statement  is valid for the discs tangent to $\cC^F$;
		
	\item for any $\eta<\eta_0$, one has that 
	\begin{itemize}
		\item[--] for any $z\in \La_\lambda$,  one considers $w_1\in f^{-1}(W_{4\eta}^{s}(f(z)))\setminus W^s_{4\eta}(z)$ and $w_2\in f(W^{s}_{4\eta}(f^{-1}(z)))$.  Let   $D_1$ and $D_2$ be two discs tangent to $\cC^F$ and  through the points $w_1, w_2$ respectively. If the diameters of $D_1$ and $D_2$ is no more than $4\eta$, then one has $D_1\cap D_2=\emptyset;$ the analogous statement is  valid for the unstable manifold of $z$;
	\item[--] for  $x\in\La$, we denote by $C_{4\eta}(x)$ the $4\eta$-tubular neighborhood of $W^{cs}_{4\eta}(x)$. Then $C_{4\eta}(x)$ is $C^1$ bi-foliated by discs tangent to the cone fields $\cC^E$ and $\cC^F$.
	
	\end{itemize}

	\end{itemize}

For $x\in\La$, we shall call  $C_{4\eta}(x)\subset B_{4\eta_0}(\Lambda)$ as a box or a cube. We denote by $\dim(E)=s$ and $\dim(F)=t$. Let   $\psi:[-1,1]^s\times [-1,1]^t\mapsto M$ be the $C^1$-embedding such that
\begin{itemize}
	\item $\psi([-1,1]^s\times[-1,1]^t)=C_{4\eta}(x)$;
	\item for $z_1\in[-1,1]^s$ and $z_2\in[-1,1]^t$, the discs $\psi(\{z_1\}\times[-1,1]^t)$ and $\psi([-1,1]^s\times \{z_2\})$ are tangent to $\cC^F$ and $\cC^E$ respectively.
	\end{itemize}
We will call $\psi(\partial([-1,1]^s)\times[-1,1]^t)$ the \emph{vertical boundary} of $C_{4\eta}(x)$ and $\psi([-1,1]^s\times\partial([-1,1]^t))$ the \emph{horizontal boundary} of $C_{4\eta}(x)$.

We say that a set  $H$ is a vertical sub-box of $C_{4\eta}(x)$, if there exists an embedding $\phi:[-1,1]^s\times [-1,1]^t\mapsto M$ such that
\begin{itemize}
	\item for any $(z_1,z_2)\in [-1,1]^s\times [-1,1]^t$, one has $\phi(\{z_1\}\times [-1,1]^t)$ and $\phi([-1,1]^s\times\{z_2\} )$ are tangent to  the cone fields $\cC^F$ and $\cC^E$ respectively.
	\item $H$ is disjoint from the vertical boundary of $C_{4\eta}(x)$;
	\item  $\phi([-1,1]^s\times \partial{([-1,1]^t)})$ is contained in the interior of the horizontal boundary of $C_{4\eta}(x)$.
\end{itemize}
One can define the horizontal sub-box analogously.

\begin{lemma}\label{l.vertical-sub-box}
	Let $x\in\La$. For $\eta<\eta_0$ and any periodic point $p\in\Lambda_\lambda\cap B_{\eta}(x)$ of period $m$,  we denote by $C(p,m)$ the connected component of $f^m(C_{4\eta}(x))\cap C_{4\eta}(x)$ containing $p.$ Then $C(p,m)$ is vertical sub-box of $C_{4\eta}(x)$. 
	\end{lemma}
\proof Since $C_{4\eta}(x)$ is foliated by discs tangent to the cone field $\cC^E$ and the discs tangent to the cone field $\cC^F$ respectively, 
$f^m(C_{4\eta}(x))$ is foliated by discs tangent to $\cC^F$; therefore $f^m(C_{4\eta}(x))\cap C_{4\eta}(x)$ is foliated by discs tangent to the cone field $\cC^E$ and the discs tangent to the cone field $\cC^F$ respectively. Notice that $f^m(C_{4\eta}(x))\cap C_{4\eta}(x)$ is saturated by discs tangent to the cone field $\cC^F$ and intersecting $f^m(W_{4\eta}^s(p))$, 
then by the choice of $\eta$,  $f^m(C_{4\eta}(x))\cap C_{4\eta}(x)$ is disjoint from the vertical boundary of $C_{4\eta}(x)$. 
Since $f^m(C_{4\eta}(x))$ is saturated by the $f^m$-images of some discs tangent to $\cC^F$, it suffices to prove the following claim.
\begin{claim}\label{c.backward-control}
	 Given  any disc $D$ intersecting  $f^m(W_{4\eta}^s(p))$ and satisfying 
	\begin{itemize}
		\item $D\subset C_{4\eta}(x)$
		\item $f^{-i}(D)$ is tangent to the cone field $\cC^F$ for $i=0,\cdots m$;
		\end{itemize}
 the diameter $f^{-m}(D)$	is less than $(\lambda+2\varepsilon)^m 4\eta$.
	\end{claim}
\proof Let $z\in D\cap f^m(W_{4\eta}^s(p))$. By the choice of $\eta$, one has that $Df^{-1}|_{T_{f^{-i}(z)}D}$ is $\varepsilon$ close to $Df^{-1}|_{F(f^{-i}(p))}$, for any $i=0,\cdots,m$.

We will inductively show that   $\diam(f^{-i}(D))<(\lambda+2\varepsilon)^i 4\eta$\: for $i=0,\cdots,m$.
It is clear that the inequality is true for $i=0$. Assume that the inequality is true for $i\leq k$. For $i=k+1$, by the choice of $\eta$, for any $y\in D$ and $j\leq k$ one has that $Df^{-1}|_{T_{f^{-j}(y)}D}$ is $\varepsilon $ close to $Df^{-1}|_{T_{f^{-j}(z)}D}$, therefore is $2\varepsilon$ close to  $Df^{-1}|_{F(f^{-j}(p))}$. Since $p\in\Lambda_\lambda$, one has  $$\diam(f^{-(k+1)}(D))\leq (\lambda+2\varepsilon)^{k+1}\diam(D)=(\lambda+2\varepsilon)^{k+1}4 \eta.$$
\endproof 
 One by-product from the proof of Lemma~\ref{l.vertical-sub-box} is the following:
\begin{corollary}\label{co:separating}
	 For each $y\in C(p,m)$, one has that $\ud(f^{-i}(y),f^{-i}(p))\leq 8\eta$, for $i=0,\cdots,m$.
	\end{corollary}
\proof For each point $y\in C(p,m)$, there exist a point $z\in W^s_{4\eta}(p)$ and a disc $D_y\subset C_{4\eta}(x)$ tangent to the cone field $\cC^F$ such that 
  $z\in D_y$ and $f^m(D_y)$ contains $y$; we denote by $\tilde{D}_y$ be the connected component of $f^m(D_y)\cap C_{4\eta}(x)$ containing $y$. One can apply Claim~\ref{c.backward-control} to $\tilde{D}_y$ and one has that $\ud(f^{-i}(y),f^{-i}(p))\leq\ud(f^{-i}(y), f^{m-i}(z))+\ud(f^{m-i}(z),f^{-i}(p))\leq 8\eta$.
\endproof

Let us continue with the proof of Theorem~\ref{thm.entropy-periodic-orbit}.
We consider the $4\eta$-cube $C_{4\eta}(p_1)$ centered at $p_1$. Since each $p_i\in B_\eta(p_1)$, by Lemma~\ref{l.vertical-sub-box},  the connected component $C(p_i,n)$ of $f^n(C_{4\eta}(p_1))\cap C_{4\eta}(p_1)$ containing $p_i$  is a vertical sub-cube of $C_{4\eta}(p_1)$.
Since the set $\{p_1,\cdots,p_k\}$ is a $(n,16\eta)$ separated set, by Corollary~\ref{co:separating}, the cubes $C(p_i,n)$ are pairwise disjoint.  By the definition of $\Lambda_\lambda$ and the choice of $\eta$, the maximal invariant set of $f^n$ in $C_{4\eta}(p_1)$ contains a hyperbolic horseshoe $K_\varepsilon$ of $k$-legs. Therefore,  
$h_{top}(f|_{K_\varepsilon})=\frac{\log k}{n}$. Now, for any point $x\in K_\e$ and any integer $l>0$, there exist $p_{i_0},\cdots, p_{i_{[l/n]-1}}$ such that $f^{jn}(x)\in C(p_{i_j},n)$ for $j=0,\cdots,[l/n]-1$.
By the choice of $\eta_0$ and Corollary~\ref{co:separating}, one has 
$$\ud\big(\frac{1}{[l/n]}\sum_{j=0}^{[l/n]-1}\frac{1}{n}\sum_{t=0}^{n-1}\delta_{f^t(f^{jn}(x))},\frac{1}{[l/n]}\sum_{j=0}^{[l/n]-1}\delta_{\cO_{p_{i_j}}}\big)<\e,$$
which implies that any accumulation $\nu$ of $\frac{1}{l}\sum_{j=0}^{l-1}\delta_{f^j(x)}$ is  $\e$-close to $\Cov(p_1,\cdots,p_k)$, then by Sigmund's result~\cite{Sig1}, all the invariant measures supported on $K_\e$ are $\e$-close to $\Cov(p_1,\cdots,p_k)$.

\endproof
\subsection{Pesin blocks and periodic points  close to non-hyperbolic ergodic measures }
In this section, we firstly define Pesin blocks close to non-hyperbolic ergodic measures, then we find periodic points close to the non-hyperbolic ergodic measures.

Let  $f$ be a partially hyperbolic diffeomorphism with one dimensional center and  $\mu$ be a non-hyperbolic ergodic measure of $f$.   We denote by $B(\mu)$ the Oseledec basin of $\mu$. For any $\chi>0$ small and any integer $k\in\mathbb{N}$, consider the following form of  Pesin block
\begin{align*}
	\Lambda_{k,\chi}=\big\{x\in B(\mu):&\:\:\frac{1}{k}e^{-n\chi}\leq\prod_{i=0}^{n-1}\norm{Df|_{E^c(f^i(x))}}\leq k e^{n\chi}, \textrm{\:for $n\in\mathbb{N}$}\\
	&\textrm{\:\:and } \ud(\frac{1}{l}\sum_{i=0}^{l-1}\delta_{f^i(x)},\mu)\leq \chi \textrm{\: for $l\geq k$}
	\big\}.
	\end{align*}

\begin{Remark}
 By definition,  one has that
 \begin{enumerate}
 	\item $\Lambda_{k,\chi}$ is a compact set (in general, $\Lambda_{k,\chi}$ is not invariant);
 \item$\Lambda_{k_2,\chi_2}\subset\Lambda_{k_1,\chi_1},\textrm{  for each $k_1\geq k_2$ and $\chi_1\geq\chi_2$.}$
    \item  $\mu(\cup_{k\in\mathbb{N}}\Lambda_{k,\chi})=1, \textrm{ for each $\chi>0$}.$
\end{enumerate}
\end{Remark}

We denote
$$H_{\chi}=\{x\in M:~\|Df^{-n}|_{E^c(x)}\|\le e^{-n\chi},\textrm{ for $n\in\N$}\}.$$

Recall that $\cU(M)$ is the set of partially hyperbolic (but roustly non-hyperbolic) diffeomorphisms with one dimensional center and two minimal strong foliations.  In the following, we will show that for an open and dense subset of $\cU(M)$, for each non-hyperbolic ergodic measure $\mu$, one can find a huge amount of hyperbolic periodic orbits in $H_\chi$ which is `close' to the measure $\mu$ and is a separated set for some size; the proof consists in combining Theorem~\ref{thm.shadow} with Liao's shadowing lemma (Lemma~\ref{l.shadowing}).
\begin{theorem}\label{p.number-of-periodic-orbit} There exists an open and dense subset $\cV(M)$ of $\cU(M)$ such that for $f\in\cV(M)$ and  a  non-hyperbolic $f$-ergodic measure $\mu$, there exists a constant $c>0$ such that 	for any $\eta$ small and any $\chi>0$, there exists   a sequence of periodic points   $\{p_1,\cdots, p_l\}\subset H_\chi$ satisfying that 
	\begin{itemize}
		\item[--] the periodic points $\{p_1,\cdots, p_l\}$ are of same period $m$ for some $m\in\mathbb{N}$;
		\item[--] $\ud(p_i,p_j)<\eta$, for any $1\leq i<j\leq l$;
		\item[--]  $\{p_1,\cdots, p_l\}$ is a $(m,16\eta)$-separated set.
		\item[--] $$    \frac{m(h_\mu(f)-\chi)}{1+c\cdot\chi}\leq\log l\textrm{\: and \:} \ud(\frac{1}{m}\sum_{j=0}^{m-1}\delta_{f^j(p_i)},\mu)<c\cdot\chi.$$
		\end{itemize}

	\end{theorem}

\proof
Let $\cV(M)$ be the open and dense subset of $\cU(M)$ given by Theorem~\ref{thm.shadow} such that for each $f\in\cV(M)$, there exists a  constant $c_0>0$  such  that for any $\varepsilon>0$, $\chi>0$   small and  any $C>1$,  there exist two integers $N(\varepsilon)\ll N(C,\chi,\varepsilon)=N $ satisfying the properties listed in Theorem~\ref{thm.shadow}.

Now, one fixes $f\in\cV(M)$ and a non-hyperbolic ergodic measure $\mu$. We fix $\chi>0$ and choose an integer $k$   such that $\mu(\Lambda_{k,\chi})>0$.
For $\eta>0$, we denote by  $K_{k,\chi}(n,\eta)$  a $(n,20\eta)$-separated set of $\Lambda_{k,\chi}$ with maximal cardinal.

By Corollary~\ref{co:entropy-separated-set}, there exists $\eta_1$ such that  for   any $\eta< \eta_1$, there exists $n$ arbitrarily large such that   $$\frac{1}{n}\log{\#K_{k,\chi}(n,\eta)}\geq h_{\mu}(f)-\frac{\chi}{2}.$$
By the uniform continuity of the center bundle, there exists $\eta_2>0$ such that for any $x,y$ with $\ud(x,y)<\eta_2$, one has 
$$e^{-\chi}<\frac{\norm{Df|_{E^c(y)}}}{\norm{Df|_{E^c(x)}}}<e^{\chi}.$$
By the compactness of $M$, there exists $\eta_3>0$ such that for $x,y$ with $\ud(x,y)<\eta_3$, one has 
$$\ud(\delta_x,\delta_y)<\chi,$$
where $\delta_x$ denotes the Dirac measure supported on $x$.
\medskip

From now on, one fixes $\eta\in(0,\min\{\eta_1,\eta_2,\eta_3 \})$.

For $e^{-4\chi}\in(0,1)$ and the splitting $TM=E^s\oplus (E^c\oplus E^u)$, by Lemma~\ref{l.shadowing}, there exist $L>0$ and $\e_0>0$ such that for any $\e\leq \e_0$, any integer $l>0$ and $x$ with the properties: 
\begin{itemize}
\item[--] $$\ud(x,f^l(x))<\e,$$
\item[--] $$\prod_{i=j}^{l-1}\norm{Df|_{E^c(f^{i}(x))}}\geq e^{4(l-k)\chi}, \textrm{  for $j=0,\cdots,l-1$;}$$
\end{itemize}
there exists a periodic point $p$ of period $l$ such that
$$\ud(f^i(p),f^i(x))<L\cdot \e, \textrm{ for $i=0,\cdots,l-1$}. $$

By Theorem~\ref{thm.shadow}, there exists a constant $c_0>0$ such that for any $\e<\min\{\e_0,\frac{\eta}{3L},\frac{\eta}{3} \}$ small, the numbers $\chi$ and $k$, there exist $N(\e)$ and $N(k,\chi,\e)$ such that for any $x\in K_{k,\chi}(n,\eta)\subset\Lambda_{k,\chi}$ and any $n+N(\e)>N(k,\chi,\e)$ (the choice of $n$ would be fixed later), there exist $z$ and $m\in\big(n+N(\e),(n+N(\e))(1+c_0\cdot \chi)\big)$ such that
\begin{itemize}
	\item[--] $$\ud(z,f^m(z))<\e;$$
	\item[--]  $$\ud(f^i(f^{N(\e)}(z)),f^i(x))<\e, \textrm{ for $i=0,\cdots,n$;}$$
	\item[--] $$\norm{Df^i|_{E^c(f^{m-i}(z))}}\geq e^{4i\chi}, \textrm{ for each $i=1,\cdots,m$.}$$
\end{itemize}
By the choice of $\e$ and Lemma~\ref{l.shadowing}, there exists a periodic point $p_x$ of period $m_x\in\big(n+N(\e),(n+N(\e))(1+c_0\cdot \chi)\big)$ such that 
$$\ud(f^i(p_x), f^i(z))<L\cdot \e<\frac{\eta}{3},\textrm{ for $i=0,\cdots,m_x-1$},$$
which implies that 
$$\ud(f^{i+N(\e)}(p_x), f^{i}(x))<\eta,\textrm{\: for $i=0,\cdots,n$}.$$
By the fact that $\eta<\eta_2$, one has 
$$\norm{Df^i|_{E^c(f^{m_x-i}(p_x))}}> e^{3i\chi}, \textrm{ for each $i=1,\cdots,m_x$.}$$

There exists an integer $N_{\e,\chi}$ such that for $n>N_\e$, one has 
$$\frac{(n+N(\e))(1+c_0\cdot \chi)-n}{n}<(c_0+1)\chi.$$

By the choice of metric on the space of probability measures and the fact that $\eta<\eta_3$, for $n\geq N_{\e,\chi}$, one has 
\begin{align*}
\ud\big(\frac{1}{m_x}\sum_{j=0}^{m_x-1}\delta_{f^j(p_x)},\mu\big)&\leq\ud(\frac{1}{m_x}\sum_{j=0}^{m_x-1}\delta_{f^j(p_x)}, \frac{1}{m_x}\sum_{j=0}^{m_x-1}\delta_{f^j(x)})+\ud(\frac{1}{m_x}\sum_{j=0}^{m_x-1}\delta_{f^j(x)},\mu) \\
&\leq \ud(\frac{1}{m_x}\sum_{j=0}^{m_x-1}\delta_{f^j(p_x)}, \frac{1}{m_x}\sum_{j=0}^{m_x-1}\delta_{f^j(x)})+\chi\\
&\leq \frac{m_x-n}{m_x}\cdot 2+\frac{n}{m_x}\cdot\chi+\chi\\
&\leq 2(c_0+2)\chi.
\end{align*}
To summarize, for $n\geq N_{\e,\chi}$ and $x\in K_{k,\chi}(n,\eta)$,  one has a hyperbolic periodic point $p_x$ of period $m_x\in\big(n+N(\e),(n+N(\e))(1+c_0\cdot \chi)\big)$ such that
\begin{itemize}
	\item $$\norm{Df^j|_{E^c(f^{m_x-j}(p_x))}}\geq e^{3j\chi}, \textrm{  for each $j=1,\cdots,m_x$;}$$
	\item  
	$$\ud(f^i(f^{N(d)}(p_x)),f^{i}(x))<\eta, \textrm{ for $i=0,\cdots,n$.}$$
	\item $$\ud\big(\frac{1}{m_x}\sum_{j=0}^{m_x-1}\delta_{f^j(p_x)},\mu\big)<2(c_0+2)\chi.$$
\end{itemize}
As a consequence, one has $p_x\in H_\chi$.

 By the compactness of $H_\chi$, there exist  $N_{\eta}$ balls $\{B(1),\cdots, B(N_\eta)\}$ of diameter $\eta$  which cover $H_\chi$.

Now, we fix an integer $n>\max\{N(k,\chi,\e),N_{\e,\chi}, N_{\eta} \}$ large such that
$$\frac{1}{n}\log{\#K_{k,\chi}(n,\eta)}\geq h_{\mu}(f)-\frac{\chi}{2}\textrm{\:\: and\:\:} \frac{1}{c_0\cdot (n+N(\e))^2\cdot \chi}>e^{-\frac{n\chi}{2}}.$$

Then we divide the set $K_{k,\chi}(n,\eta)$ into at most $ c_0\cdot (n+N(\e))\cdot\chi$ subsets such that in each of these subsets, $m_x$ is constant. Let  $\tilde{K}_{k,\chi}(n,\eta)$ be one of such  subsets of $K_{k,\chi}(n,\eta)$ with the maximal cardinal, then one has
$$\#\tilde{K}_{k,\chi}(n,\eta)\geq \frac{1}{c_0\cdot (n+N(\e))\cdot \chi}\#K_{k,\chi}(n,\eta)\geq  \frac{1}{c_0\cdot (n+N(\e))\cdot \chi}\cdot  e^{n\cdot(h_{\mu}(f)-\frac{\chi}{2})}.$$

\begin{claim} For the points  $x,y\in\tilde{K}_{k,\chi}(n,\eta)$, one has that $p_x=p_y$ if and only if $x=y$.
\end{claim}
\proof If one has  $p_x=p_y$ for some $x,y\in \tilde{K}_{k,\chi}(n,\eta)$,   one has that 
$$\ud(f^i(x),f^i(y))\leq \ud(f^i(x),f^{N(\e)+i}(p_x))+\ud(f^i(y),f^{N(\e)+i}(p_y))<2\eta \textrm{ for $i=0,\cdots,n$,}$$  which contradicts to the choice of $x,y$.
\endproof
Since $\tilde{K}_{k,\chi}(n,\eta)$ is also a $(n,20\eta)$-separated set, for different points  $x,y\in\tilde{K}_{k,\chi}(n,\eta)$, one has that there exists $j\in[N(\e),n+N(\e)]$ such that
$$\ud(f^j(p_x),f^j(p_y))>16\eta.$$

Now, there exists $j\in\{1,2,\cdots, N_\eta\}$ such that $B(j)\cap\{p_x:x\in\tilde{K}_{k,\chi}(n,\eta)\}$ has the maximal cardinal, then one has
\begin{align*} \#(B(j)\cap\{p_x:x\in\tilde{K}_{k,\chi}(n,\eta)\})&\geq \frac{1}{N_{\eta}} \frac{1}{c_0\cdot (n+N(\e))\cdot \chi}\cdot e^{n\cdot(h_{\mu}(f)-\frac{\chi}{2})}\\
&>\frac{1}{c_0\cdot (n+N(\e))^2\cdot \chi}\cdot e^{n\cdot(h_{\mu}(f)-\frac{\chi}{2})}\\
&>e^{n\cdot(h_{\mu}(f)-\chi)}.
\end{align*}
We denote $m$ as the period of $p_x$ for $x\in\tilde{K}_{k,\chi}(n,\eta),$ then by the choice of $n$, one has that 
$$\log \#(B(j)\cap\{p_x:x\in\tilde{K}_{k,\chi}(n,\eta)\})\geq \frac{m\cdot(h_{\mu}(f)-\chi)}{1+(c_0+1)\chi}.$$

  Now we only need to set $c=2(c_0+2)$ and denote $B(j)\cap\{p_x:x\in\tilde{K}_{k,\chi}(n,\eta)\}=\{p_1,\cdots,p_l\}$, then  $\{p_1,\cdots,p_l\}$ is the set of periodic points with the posited properties.
\endproof

\section{Non-hyperbolic ergodic measures approached by hyperbolic basic sets}

In this section, using the results from previous sections, we give the proof of our main results.
\proof[Proof of Theorem~\ref{thmA}]
Let $\cV(M)$ be the open and dense subset of $\cU(M)$ given by Theorem~\ref{p.number-of-periodic-orbit}. One fixes  $f\in\cV(M)$ and an $f$-ergodic measure $\mu$.

If $\mu$ is hyperbolic, by Katok-Gelfert theorem ~\cite[Theorem 1]{Ge}, $\mu$ is approached by hyperbolic sets in weak$*$-topology and in entropy. 

Now, we consider the case that $\mu$ is non-hyperbolic, that is, the center Lyapunov exponent of $\mu$ vanishes. 

For $\chi>0$, the splitting $TM=E^s\oplus (E^c\oplus E^u)$ and $\lambda=e^{-\chi}\in(0,1)$, by Theorem~\ref{thm.entropy-periodic-orbit},  there exists $\eta_0=\eta_0(\chi)>0$ such that for any $\eta<\eta_0$, any $n\in\mathbb{N}$ and  any periodic points $\{q_1,\cdots, q_k \}\subset H_\chi$ of   period $n$, if one has 
\begin{itemize}
	\item $\ud(q_i,q_j)<\eta$;
	\item $\{q_1,\cdots, q_k \}$ is a $(n,16\eta)$-separated set;
\end{itemize} 
then there exists a hyperbolic set $K_\chi$ of index $\dim(E^s)$ such that 
\begin{itemize}
	\item[--] 
	$h_{top}(f|_{K_\chi})\geq \frac{\log k}{n}$;
	\item[--] for any  $\nu\in\cM_{inv}(K_\chi,f)$, one has  $\ud(\nu,\Cov(q_1,\cdots,q_k))<\chi$.
\end{itemize}

By Theorem~\ref{p.number-of-periodic-orbit},  for the non-hyperbolic measure $\mu$, there exists a constant $c>0$ such that 	for   $\eta<\eta_0$ small and  the number $\chi$, there exists a sequence of periodic points   $\{p_1,\cdots, p_l\}\subset H_\chi$ satisfying that 
\begin{itemize}
	\item[--] periodic points   $\{p_1,\cdots, p_l\}$ are of same period $m$ for some $m\in\mathbb{N}$;
	\item[--] $\ud(p_i,p_j)<\eta$, for any $1\leq i<j\leq l$;
	\item[--] $\{p_1,\cdots, p_l\}$ is a $(m,16\eta)$-separated set.
	\item[--] $$    \frac{m(h_\mu(f)-\chi)}{1+c\cdot\chi}\leq \log l\textrm{\: and \:} \ud(\frac{1}{m}\sum_{j=0}^{m-1}\delta_{f^j(p_i)},\mu)<c\cdot\chi.$$
\end{itemize}
Then Theorem~\ref{thm.entropy-periodic-orbit} provides a hyperbolic set $K_\chi$ such that 
\begin{itemize}
	\item $h_{top}(f|_{K_\chi})\geq \frac{h_\mu(f)-\chi}{1+c\cdot\chi}$;
	\item for  $\nu\in\cM_{inv}(K_\chi,f)$, one has 
	$\ud(\nu,\Cov(p_1,\cdots,p_l))<\chi$, which implies $\ud(\nu,\mu)<(c+1)\chi.$
\end{itemize}

 By Theorem~\ref{thm.middle-value-entropy}, one gets a hyperbolic set $\tilde{K}_\chi\subset K_\chi$ which is  $(c+1)\cdot\chi$-close to $\mu$ in weak$*$-topology and in entropy. By the arbitrariness of $\chi$, one has that $\mu$ is approached by hyperbolic sets of index $\dim(E^s)$ in weak$*$-topology and in entropy. 
 
 Applying the analogous argument for the reversed dynamics in $\cU(M)$, one has that for an open and dense subset in $\cU(M)$, each non-hyperbolic ergodic measure is also appraoched by hyperbolic sets of index $\dim(E^s)+1$ in weak$*$-topology and in entropy.

\endproof

Now, using Theorem~\ref{thmA}, we  give the proof of Corollary~\ref{co:path-entropy}.
\proof[Proof of Corollary~\ref{co:path-entropy}]
One only needs to show that given a hyperbolic ergodic measure $\mu$ and a non-hyperbolic ergodic measure $\nu$, there exists a path connecting them. 
By Theorem~\ref{thmA}, there exist sequences of periodic orbits $\{p_n\}_{n\geq 1}$ and $\{q_n\}_{n\geq 1}$ of same index such that 
$\delta_{\cO_{p_n}}$ converges to $\mu$ and $\delta_{\cO_{q_n}}$ converges to $\nu$. By the minimality of strong foliations, the periodic points $p_n$ and $q_n$ are pairwise homoclinically related. Since $\cO_{q_{n+1}}$ and $\cO_{q_n}$ are homoclinically related, by Smale-Birkhoff theorem, there exists a hyperbolic horseshoe $\La_n$ containing $q_{n+1}$ and $q_n$, then by Theorem B in ~\cite{Sig2} (see also the comments after it), there exists a path $\{\nu_t\}_{t\in[1-3^{-n}, 1-3^{n+1}]}$ in the set of ergodic measures supported on $\La_n$ which connects $\delta_{\cO_{q_n}}$ to $\delta_{\cO_{q_{n+1}}}$. Analogously, in the set of ergodic measures,  one has paths $\{\nu_t\}_{t\in [3^{-n-1},3^{-n}]}$ which connects  $\delta_{\cO_{p_{n+1}}}$ to $\delta_{\cO_{p_{n}}}$ and $\nu_{t\in[\frac{1}{3},\frac{2}{3}]}$ which connects  $\delta_{\cO_{p_1}}$ to $\delta_{\cO_{q_{1}}}$. Let $\nu_0=\mu$ and $\nu_1=\nu$, then $\{\nu_t\}_{t\in[0,1]}$ gives a path connecting $\mu$ to $\nu$. 

Since on each hyperbolic sets, there always exist measures of maxiaml entropy (see for instance ~\cite{Bo}), by Theorem~\ref{thmA}, the set of hyperbolic ergodic measures is entropy dense in the set of ergodic measures.

\endproof

\bibliographystyle{plain}

\vskip 5pt

\begin{tabular}{l l l}
	\emph{\normalsize Dawei Yang}
	& \quad\quad  &
	\emph{\normalsize Jinhua Zhang}
	\medskip\\
	
	\small School of Mathematical Sciences
	&& \small Laboratoire de Math\'ematiques d'Orsay\\
	\small Soochow University
	&& \small CNRS - Universit\'e Paris-Sud\\
	\small Suzhou, 215006, P.R. China
	&& \small  Orsay 91405, France\\
	\small \texttt{yangdw1981@gmail.com, yangdw@suda.edu.cn}
	&& \small \texttt{Jinhua.zhang@math.u-psud.fr}
\end{tabular}

\end{document}